\documentclass[11pt,a4paper]{article}

\usepackage[paperwidth=22cm]{geometry}
\usepackage[T1]{fontenc}
\usepackage[latin1]{inputenc}
\usepackage{graphicx}
\usepackage{amsmath}
\usepackage{amssymb}

\newtheorem{theorem}{Theorem}
\newtheorem{corollary}{Corollary}

\newcommand{\bl}{\boldsymbol}

\newcommand{\lef}{\left(}
\newcommand{\rig}{\right)}
\newcommand{\ph}{\phantom}
\newcommand{\la}{\langle}
\newcommand{\ra}{\rangle}
\newcommand{\ip}{\lrcorner}
 \newcommand{\qed}{\hfill \ensuremath{\Box}}

\begin{document}

\title{\textbf{Pure Subspaces, Generalizing the Concept \\ of Pure Spinors}}

\author{\textbf{Carlos  Batista}\\
\small{Departamento de F\'{\i}sica}\\
\small{Universidade Federal de Pernambuco}\\
\small{50670-901 Recife-PE, Brazil}\\
\small{carlosbatistas@df.ufpe.br}}
\date{}

%\date{\small{\today}}

%Recife - PE, Brazil}
%\email{carlosbatistas@df.ufpe.br}
%\affiliation{Departamento de Física, Universidade Federal de Pernambuco, 50670-901
%Recife - PE, Brazil}

%\date{\today}

%\thanks{CNPq financial support.}
%\archivePrefix{}
%\keywords{Isotropic Spaces, Pure Spinors, Integrability, Clifford Algebra, Twistors}

\maketitle
\begin{abstract}
The concept of pure spinor is generalized, giving rise to the notion of pure subspaces, spinorial subspaces associated to isotropic vector subspaces of non-maximal dimension. Several algebraic identities concerning the pure subspaces are proved here, as well as some differential results. Furthermore, the freedom in the choice of a spinorial connection is exploited in order to relate twistor equation to the integrability of maximally isotropic distributions. \textsl{(Keywords: Isotropic Spaces, Pure Spinors, Integrability, Clifford Algebra, Twistors)}
\end{abstract}

\section{Introduction}

It is well-known that given a spinor $\hat{\varphi}$ one can construct a vector subspace $N_{\hat{\varphi}}$ spanned by the vectors that annihilate $\hat{\varphi}$ under the Clifford action, $v\in N_{\hat{\varphi}}$ whenever $v\cdot\hat{\varphi}= 0$. This subspace is necessarily isotropic, which means that $\la v,v\ra= 0$ for all $v\in N_{\hat{\varphi}}$. Particularly, if the dimension of $N_{\hat{\varphi}}$ is maximal $\hat{\varphi}$ is said to be a pure spinor. Thus, to every spinor it is associated an isotropic vector subspace. The idea of the present article is to go the other way around and associate to every isotropic subspace $I$ a spinorial subspace $\hat{L}_I$. As we shall see, these spinorial subspaces provide a natural generalization for the concept of pure spinors and, therefore, we shall say that $\hat{L}_I$ is the pure subspace associated to $I$. As a consequence, some classical results are shown to be just particular cases of broader theorems. Hence, this study can shed more light on the role played by the pure spinors on physics and mathematics.

Spinors have been around for a long time, they were first discovered a century ago by \'{E}lie Cartan \cite{Cartan0}. Since then, they led to great achievements in a multitude of areas of mathematics and physics. For instance, spinors are of fundamental importance in particle physics, since fermions are represented by spinorial fields (see \cite{Fermions-Spinors} for an alternative treatment). In addition, spinors can often be used as helpful tools for calculating physical and geometrical quantities. As examples let us recall Witten's elegant proof of energy positivity in general relativity \cite{Witten-Energy} and the calculation of scattering amplitudes in quantum field theory \cite{scattering4d-1}. From the geometrical point of view, a particularly important class of spinors is formed by the so-called pure spinors. Such spinors are also acquiring increasing significance in high-energy physics, as illustrates the pure spinor formalism in string theory \cite{Nathan}. In spite of such undeniable usefulness, the role played by spinors is certainly not completely understood yet. Hopefully, the work presented here will contribute for a further understanding of spinorial calculus and, specially, of the meaning of pure spinors.

The isotropic subspaces are of great relevance in mathematical-physics, a fact that became clear after the works by R. Penrose wherein he introduced the spinorial calculus and the null tetrad formalism in 4-dimensional general relativity \cite{Penrose-Null}. Such approaches, in which isotropic distributions play a prominent role, brought much progress to general relativity. Notably, Kinnersley was able to find analytically all solutions of Einstein's vacuum equation in 4 dimensions for the case of the space-time admitting two independent integrable distributions of isotropic planes \cite{typeD,Plebanski2,art2}. Particularly, all known 4-dimensional black-holes are contained in this class of solutions. Higher-dimensional manifolds have been a subject of great importance in theoretical physics, not to mention in mathematics. Following the successful track of 4-dimensional general relativity recent works have shown that isotropic structures do also play an important role in higher dimensions. In reference \cite{TaChabert-KY} it was proved a relation between integrable isotropic structures and the existence of a conformal Killing-Yano tensor, while in \cite{HigherGSisotropic} these structures were used to classify the Weyl tensor as well as to partially generalize the Goldberg-Sachs theorem, see also \cite{art4}.

The outline of this article is as follows. In Section \ref{section V+V*} it is introduced the notation adopted throughout the article to deal with spinors. Section \ref{section-AlgebRes} defines the pure subspaces and present a theorem on algebraic identities satisfied by them. Then, Section \ref{section-connection} extends the formalism from vector spaces to fiber bundles over differential manifolds. In addition, this section reviews the issue of introducing a connection on the spinorial bundle and explore the freedom in its choice. Section \ref{section-Integrability} presents some differential results regarding the pure subspaces, connecting the integrability of an isotropic distribution with differential constraints satisfied by the elements of the associated pure subspace. Finally, in Section \ref{section-twistor} it is proved a theorem relating twistors to the existence of integrable maximally isotropic distributions. All over this work it is being assumed that the vector spaces and the manifolds are even-dimensional and endowed with a non-degenerate metric of arbitrary signature. Unless otherwise stated, Einstein summation convention is used, meaning that repeated indices are summed. The results are all local.

%%%%%%%%%%%%%%%%%%%%%%%%%%%%%%%%%%%%%%%%%%%%%%%%%%%%%%%%%%%%%%%%%%%%%%%%%%%%%%%%%%%%%%%%%%%%%%%%%%%%%
%%%%%%%%%%%%%%%%%%%%%%%%%%%%%%%%%%%%%%%%%%%%%%%%%%%%%%%%%%%%%%%%%%%%%%%%%%%%%%%%%%%%%%%%%%%%%%%%%%%%%
%%%%%%%%%%%%%%%%%%%%%%%%%%%%%%%%%%%%%%%%%%%%%%%%%%%%%%%%%%%%%%%%%%%%%%%%%%%%%%%%%%%%%%%%%%%%%%%%%%%%%
%%%%%%%%%%%%%%%%%%%%%%%%%%%%%%%%%%%%%%%%%%%%%%%%%%%%%%%%%%%%%%%%%%%%%%%%%%%%%%%%%%%%%%%%%%%%%%%%%%%%%
%%%%%%%%%%%%%%%%%%%%%%%%%%%%%%%%%%%%%%%%%%%%%%%%%%%%%%%%%%%%%%%%%%%%%%%%%%%%%%%%%%%%%%%%%%%%%%%%%%%%%
%%%%%%%%%%%%%%%%%%%%%%%%%%%%%%%%%%%%%%%%%%%%%%%%%%%%%%%%%%%%%%%%%%%%%%%%%%%%%%%%%%%%%%%%%%%%%%%%%%%%%
%%%%%%%%%%%%%%%%%%%%%%%%%%%%%%%%%%%%%%%%%%%%%%%%%%%%%%%%%%%%%%%%%%%%%%%%%%%%%%%%%%%%%%%%%%%%%%%%%%%%%
%%%%%%%%%%%%%%%%%%%%%%%%%%%%%%%%%%%%%%%%%%%%%%%%%%%%%%%%%%%%%%%%%%%%%%%%%%%%%%%%%%%%%%%%%%%%%%%%%%%%%

\section{Spinors and the Space $V\oplus V^*$} \label{section V+V*}

Given a vector space $\mathcal{V}$ endowed with a non-degenerate inner product $\la\,,\ra$, the Clifford Algebra $Cl(\mathcal{V})$ is an algebra in this vector space such that:
\begin{equation}\label{Clifford Product}
    v\, u \,+\, u\, v \,=\, 2\, \la v, u \ra \,\;\;\; \forall \; v,u\,\in\;\mathcal{V}.
\end{equation}
The space of spinors associated to $(\,\mathcal{V} , \la\,,\ra\,)$ is a vectorial space $S$ where an irreducible and faithful representation of $Cl(\mathcal{V})$ acts. In even dimensions it is always possible to find a matrix representation for $Cl(\mathcal{V})$, if $\dim(\mathcal{V})=2n$ then the least-dimensional faithful representation of this algebra is provided by $2^n\times2^n$ matrices. Therefore, in this case spinors are represented by the column vectors on which these matrices act, so that the space of spinors has dimension $2^n$.

Now, let us deal with the vector space $V\oplus V^*$, where $V$ is an $n$-dimensional vector space and $V^*$ is  its dual. The space $V\oplus V^*$ is naturally endowed with a non-degenerate symmetric inner product $\la\,,\ra$ defined by:
\begin{equation}\label{Inner Product}
    \la e+\theta, e'+ \theta'\ra \,\equiv\, \frac{1}{2}\left[\,\theta(e') \,+\, \theta'(e)\,\right]\;\;\;\;\;\; e,e'\in V \;\textrm{ and }\;  \theta, \theta'\in V^*\,.
\end{equation}
Let us introduce the basis $\{e_1,e_2,...,e_n\}$ for $V$ and denote by $\{\theta^1,\theta^2,...,\theta^n\}$ the dual basis, $\theta^i(e_j)=\delta_j^i$. So, thanks to (\ref{Inner Product}) we have that the following relations hold:
 $$  \la e_i,e_j\ra  \,=\, 0\,, \quad\;\quad  \la e_i,\theta^j\ra \,=\, \frac{1}{2}\,\delta_i^j\,,  \quad\;\quad  \la\theta^i, \theta^j\ra \,=\,0\,.$$
The interesting feature of the space $V\oplus V^*$ endowed with this natural inner product is that the spinors can be constructed quite easily. Indeed, the algebra $Cl(V\oplus V^*)$ admits a representation in the space $\wedge V^*$, the exterior algebra of $V^*$. The action of $V\oplus V^*$ vectors in this representation being defined by:
\begin{equation}\label{Representation}
    (e+\theta)\cdot \varphi \,\equiv\, e\ip\varphi \,+\, \theta \wedge \varphi  \quad\quad\; e\in V\; , \; \theta \in V^* \;, \; \varphi\in \wedge V^*\,,
\end{equation}
where $e\ip\varphi$ means the interior product of the vector $e$ on the form $\varphi$. The action of a scalar is defined in the trivial way, $\lambda\cdot\varphi = \lambda\varphi$. In order to see that this is indeed a representation of the Clifford algebra, note that by successive application of some vector $v= e+\theta$ on the spinor $\varphi$ we get
\begin{align*}
 v\cdot\left[v\cdot\varphi\right] \,&=\, (e+\theta)\cdot\left[e\ip\varphi \,+\, \theta \wedge \varphi\right] \\
 &=\,e\ip e\ip\varphi \,+\,  e\ip( \theta \wedge \varphi) \,+\, \theta \wedge (e\ip\varphi) \,+\, \theta \wedge\theta \wedge \varphi\\
  &=\,\theta(e)\, \varphi\,=\,  \la v,v\ra \,\varphi \,=\,  (v\,v)\cdot\varphi \,,
\end{align*}
where in the last step it was used Eq. (\ref{Clifford Product}). Since the space $\wedge V^*$ has $2^n$ dimensions we conclude that this is, indeed, the space of spinors, $S=\wedge V^*$ \cite{Gualtieri}.

To simplify the notation, it will be used the abbreviation $\theta^{12...k} \equiv \theta^1\wedge\theta^2\wedge...\wedge\theta^k$. Furthermore, in order to avoid  any confusion it is worth distinguishing the elements of $Cl(V{\oplus}V^*)$  from the spinors of $S=\wedge V^*$. With this intent from now on we shall always denote the spinors by a Greek letter with a hat on top. For instance, $\theta^i$ is an element of $V\oplus V^*$ while $\hat{\theta}^i$ is a spinor.

In the present formalism an inner product can be easily introduced on the space of spinors. Such product is non-degenerate and defined up to an arbitrary scale factor. This scale can be fixed by an arbitrary choice of basis for $S$, here let us choose the basis $\{\hat{1},\hat{\theta}^i,\hat{\theta}^{ij}, ..., \hat{\theta}^{12...n}\}$. Then the inner product $( \,,):S\times S\rightarrow \mathbb{C}$ is defined to be such that:
\begin{equation}\label{InnerSpinors}
    < \,\hat{\varphi}^{\,t} \wedge\hat{\psi}\, >_n \,\,=\, \lambda\, \hat{\theta}^{12...n}\;\;\; \Longleftrightarrow\;\;\; (\hat{\varphi},\hat{\psi})\,\equiv\, \lambda\,.
\end{equation}
Where $<\hat{\varphi}>_k$ means the component of degree $k$ of the form $\hat{\varphi}$, while $\hat{\varphi}^{\,t}$ means the reverse of $\hat{\varphi}$. For instance, if $\hat{\varphi}= \hat{1}+ \hat{\theta}^{12...k}$ then $<\hat{\varphi}>_k= \hat{\theta}^{12...k}$ and $\hat{\varphi}^{\,t}= \hat{1}+ \hat{\theta}^{k...21}$. It is not difficult convincing ourselves that $(v\cdot\hat{\varphi},\hat{\psi})= (\hat{\varphi},v\cdot\hat{\psi})$ for any $v\in V\oplus V^*$. Such property implies that this inner product is invariant under the group $Spin_+(V\oplus V^*)$, the double cover of $SO_+(V\oplus V^*)$. More generally, if $\Omega\in Cl(V\oplus V^*)$ then it can be proved that $(\Omega\cdot\hat{\varphi},\hat{\psi})= (\hat{\varphi},\Omega^t\cdot\hat{\psi})$.

If $v_1,\, \ldots,v_p$ are vectors in $V\oplus V^*$ then we shall define $v_1\wedge\ldots\wedge v_p$ to be the complete antisymmetric part of the Clifford product $v_1v_2\cdots v_p$. In particular, we have that $v_1\wedge v_2= \frac{1}{2!}(v_1v_2- v_2v_1)$. Then, the so-called pseudo-scalar of the Clifford algebra $Cl(V\oplus V^*)$ is defined by:
$$  \mathcal{I} \,=\, 2^n\,e_1\wedge\theta^1\wedge e_2\wedge\theta^2\wedge \ldots \wedge e_n\wedge\theta^n\,.  $$
It is simple matter to prove that $ \mathcal{I}\, \mathcal{I}=1$. So that such object can be used to split the spinor space as $S=S^+\oplus S^-$, where $\hat{\varphi}\,\in\,S^{\pm}$ if $ \mathcal{I}\cdot\hat{\varphi}=\pm\hat{\varphi}$. The elements of $S^+$ are called the Weyl spinors of positive chirality, while those of $S^-$ are the Weyl spinors of negative chirality. It can be shown that $S^+$ is spanned by the spinors of even degree, while $S^-$ is spanned by the spinors of odd degree, so that $\dim(S^\pm)= 2^{n-1}$. More about the spinors of $Cl(V\oplus V^*)$ can be found in \cite{Gualtieri}, where this formalism is applied to the space $TM\oplus TM^*$, with $TM$ being the tangent bundle of a manifold $M$. In order for the reader to get acquainted with the language introduced so far, let us work out a simple example.
\\
\\
\small{
\textit{\textbf{Example 1}}\\
When $\dim(V)=2$ we have that $\{e_1,e_2,\theta^1,\theta^2\}$ provides a basis for $V\oplus V^*$. Thus, a basis for the Clifford algebra $Cl(V\oplus V^*)$ is formed by the following 16 elements:
$$\{1,\,e_i,\,\theta^i,\, e_1\wedge e_2,\, e_i\wedge \theta^j,\, \theta^1\wedge\theta^2,\, e_i\wedge\theta^1\wedge\theta^2,\, e_1\wedge e_2\wedge \theta^i,\, e_1\wedge e_2\wedge\theta^1\wedge\theta^2\}\,.$$
While a basis for the spinor space is given by $\{\hat{1},\hat{\theta}^1, \hat{\theta}^2, \hat{\theta}^{12}\}$. Since $\hat{1}$ and $\hat{\theta}^{12}$ are forms of even degree they are Weyl spinors of positive chirality, while $\hat{\theta}^1$ and $\hat{\theta}^2$ have negative chirality. Note, for instance, that the following relations hold:
$$ e_i\cdot\hat{1} = 0\,, \quad\; e_i\cdot\hat{\theta}^j\,=\, \delta^j_i\,\hat{1} \,, \quad\; e_1\cdot\hat{\theta}^{12}\,=\,\hat{\theta}^{2} \,, \quad\; \theta^i\cdot\hat{1} = \hat{\theta}^i \,, \quad\;  \theta^1\cdot\hat{\theta}^2= \hat{\theta}^{12} \,, \quad\;  \theta^1\cdot\hat{\theta}^{12}= 0\,.$$
Moreover, the non-zero inner products are given by:
$$  (\hat{1},\hat{\theta}^{12}) \,=\, -(\hat{\theta}^{12},\hat{1})\,=\, 1  \quad \textrm{and} \quad
(\hat{\theta}^1, \hat{\theta}^2) \,=\, -(\hat{\theta}^2, \hat{\theta}^1) \,=\, 1\,.   $$
In particular, note that this inner product is skew-symmetric.\qed}
\normalsize\\

Although it may appear too restrictive working with vector spaces of the form $V{\oplus}V^*$, this is not the case at all. Every even-dimensional vector space endowed  with a non-degenerate metric can easily be cast in the form $V\oplus V^*$ when complexified. For example, in the Minkowski space $\mathbb{R}^{1,3}$ a complex null tetrad basis can always be introduced. In the standard notation of General Relativity this null tetrad is denoted by $\{l,m,n,\overline{m}\}$, with $\la l,n\ra =1$ and $\la m,\overline{m}\ra= -1$, all other inner products between the basis vectors being zero. So it is possible to make the following associations: $l\leftrightarrow e_1$, $m\leftrightarrow e_2$, $n\leftrightarrow2\theta^1$, $\overline{m}\leftrightarrow-2\theta^2$. Thus, $\{l,m\}$ can be seen as a basis for $V$ and  $\{n,\overline{m}\}$ a basis for $V^*$. So, from now on all calculations will be done on vector spaces of the form $V\oplus V^*$ and over the complex field. But it must be clear that the results can be easily carried to all complexified even-dimensional spaces. When the metric on the even-dimensional space has split signature the results are also valid without complexification, since in this case the isotropic subspaces can have dimension equal to half of the dimension of the full vector space. The results on real even-dimensional vector spaces can be extracted from the complex case by choosing suitable reality conditions, in the spirit of \cite{art2,art4}.

%%%%%%%%%%%%%%%%%%%%%%%%%%%%%%%%%%%%%%%%%%%%%%%%%%%%%%%%%%%%%%%%%%%%%%%%%%%%%%%%%%%%%%%%%%%%%%%%%%%%%
%%%%%%%%%%%%%%%%%%%%%%%%%%%%%%%%%%%%%%%%%%%%%%%%%%%%%%%%%%%%%%%%%%%%%%%%%%%%%%%%%%%%%%%%%%%%%%%%%%%%%
%%%%%%%%%%%%%%%%%%%%%%%%%%%%%%%%%%%%%%%%%%%%%%%%%%%%%%%%%%%%%%%%%%%%%%%%%%%%%%%%%%%%%%%%%%%%%%%%%%%%%
%%%%%%%%%%%%%%%%%%%%%%%%%%%%%%%%%%%%%%%%%%%%%%%%%%%%%%%%%%%%%%%%%%%%%%%%%%%%%%%%%%%%%%%%%%%%%%%%%%%%%
%%%%%%%%%%%%%%%%%%%%%%%%%%%%%%%%%%%%%%%%%%%%%%%%%%%%%%%%%%%%%%%%%%%%%%%%%%%%%%%%%%%%%%%%%%%%%%%%%%%%%
%%%%%%%%%%%%%%%%%%%%%%%%%%%%%%%%%%%%%%%%%%%%%%%%%%%%%%%%%%%%%%%%%%%%%%%%%%%%%%%%%%%%%%%%%%%%%%%%%%%%%
%%%%%%%%%%%%%%%%%%%%%%%%%%%%%%%%%%%%%%%%%%%%%%%%%%%%%%%%%%%%%%%%%%%%%%%%%%%%%%%%%%%%%%%%%%%%%%%%%%%%%
%%%%%%%%%%%%%%%%%%%%%%%%%%%%%%%%%%%%%%%%%%%%%%%%%%%%%%%%%%%%%%%%%%%%%%%%%%%%%%%%%%%%%%%%%%%%%%%%%%%%%

\section{Pure Spinorial Subspaces, Algebraic Results }\label{section-AlgebRes}

While dealing with spinors it is common, and often valuable, associating to every non-zero spinor $\hat{\varphi}$ a vector subspace
$N_{\hat{\varphi}}\subset (V\oplus V^*)$  defined by:
$$ N_{\hat{\varphi}} = \{\, v \in (V\oplus V^*)\;|\; v\cdot \hat{\varphi} = 0 \,\}\,.$$
These vector subspaces are called isotropic or totally null, because every vector belonging to $N_{\hat{\varphi}}$ has zero norm. This can be easily verified, if $v\in N_{\hat{\varphi}}$ then $\la v,v\ra\hat{\varphi} = (vv)\cdot\hat{\varphi}= v\cdot (v\cdot\hat{\varphi})= 0$, so that $\la v,v\ra= 0$. The maximum dimension that an isotropic subspace of $V\oplus V^*$ can have is $n=\dim(V)$. Therefore, an isotropic subspace of dimension $n$ is called a maximally isotropic subspace.  In the particular case of $N_{\hat{\varphi}}$ being maximally isotropic, the spinor $\hat{\varphi}$ is said to be a pure spinor. For example, in the notation of last section the spinor $\hat{1}$ is such that $N_{\hat{1}}= \textrm{Span}\{e_1,e_2,...,e_n\}$, so that $\hat{1}$ is a pure spinor. The pure spinors are very special objects in mathematics and have been studied since the beginning of XX century, more about them can be found in \cite{Cartan,Chevalley,Spinor-thesis}. Pure spinors are also acquiring increasing relevance in physics, particularly in string theory \cite{Nathan}.

However, it seems to have been overlooked that it is also possible to go the other way around and associate spinors to isotropic subspaces. More precisely, given an isotropic subspace $I$, one can define a spinorial subspace $\hat{L}_{I}$ spanned by all spinors annihilated by the Clifford action of $I$. The aim of the present article is to explore this new path.

Suppose that $I\subset (V\oplus V^*)$ is an isotropic subspace. Then, let us define the subspace $\hat{L}_{I}\subset S$ as follows:
\begin{equation}\label{Spinor subspace}
    \hat{L}_{I} \,\equiv\, \left\{\,\hat{\varphi} \in S \;|\; v\cdot \hat{\varphi}= 0 \;\; \forall \;v\in I\, \right\} \,.
\end{equation}
It is trivial to see that $\hat{L}_{I}$ is, indeed, a vector subspace of the spinor space $S=\wedge V^*$. In the particular case of $I$ being maximal, $\dim(I)=n$, it follows that $\hat{L}_{I}= \textrm{Span}\{\hat{\psi}\}$, where $\hat{\psi}$ is the pure spinor associated to $I$ and $\textrm{Span}\{\hat{\psi}\}$ is the one-dimensional subspace spanned by it. This section is devoted to enunciate and prove several algebraic results concerning the spaces $\hat{L}_{I}$, where $I$ is any totaly null subspace of $V \oplus V^*$. Before proceeding let us see a simple example.
\\
\\
\small{
    \textit{\textbf{Example 2}}\\
    If $\dim(V)=3$ we have that $\{e_1,e_2,e_3,\theta^1,\theta^2,\theta^3,\}$ is a basis to $V \oplus V^*$. The spinor space is generated by the basis $\{\hat{1},\hat{\theta}^1, \hat{\theta}^2, \hat{\theta}^3, \hat{\theta}^{12},\hat{\theta}^{13},\hat{\theta}^{23},\hat{\theta}^{123} \}$. Then, let us define the following isotropic subspaces of $V \oplus V^*$:
    $$ I_1 \,=\,\textrm{Span}\{e_1\}\,,\quad I_2\,=\, \textrm{Span}\{e_1,\theta^2\}\,,\quad  I_3\,=\, \textrm{Span}\{e_1,\theta^2, \theta^3\}\,. $$
    Thus, the associated spinorial subspaces are respectively given by:
    $$ \hat{L}_{I_1}\,=\, \textrm{Span}\{\hat{1},\hat{\theta}^2, \hat{\theta}^3,\hat{\theta}^{23} \} \,,\quad  \hat{L}_{I_2}\,=\, \textrm{Span}\{\hat{\theta}^2,\hat{\theta}^{23}\} \,,\quad  \hat{L}_{I_3}\,=\, \textrm{Span}\{\hat{\theta}^{23} \}\,. $$\qed
   }\normalsize
\\
\\
Now, let us enunciate, in the form of a theorem, several interesting algebraic results regarding the spinorial subspaces $\hat{L}_{I}$.
\begin{theorem}\label{Theor-Algebraic}
Given two isotropic subspaces $I,I' \subset (V \oplus V^*)$, with $V$ being a complexified vector space of dimension $n$, then it follows that:\\
\textbf{(1)} \quad $\dim(\hat{L}_{I}) \,=\, 2^{n-\dim(I)}$\\
\textbf{(2)}\quad $\hat{L}_{I} = \hat{L}_{I'}$ $\Longleftrightarrow$ $I=I'$ \\
\textbf{(3)}\quad $(\hat{L}_{I} + \hat{L}_{I'}) \subset \hat{L}_{I\cap I'}$\\
\textbf{(4)}\quad If $I''=I + I'$ is isotropic then $\hat{L}_{I+I'} = \hat{L}_{I} \cap\hat{L}_{I'}$\\
\textbf{(5)}\quad If $u\in(V \oplus V^*)$ is such that $u\cdot \hat{\varphi}=0 $ $\,\forall\,\;\hat{\varphi} \in \hat{L}_{I}$ then $u\in I$\\
\textbf{(6)}\quad $I'\subset I$ $\Longleftrightarrow$ $\hat{L}_{I} \subset \hat{L}_{I'}$\\
\textbf{(7)}\quad If $I\neq\{0\}$ then  $(\hat{\varphi},\hat{\psi})=0$ $\,\forall\,\;\hat{\varphi},\hat{\psi} \in \hat{L}_{I}$ \\
\textbf{(8)}\quad $I''=I+I'$ is isotropic $\Longleftrightarrow$  $\hat{L}_{I} \cap\hat{L}_{I'}\neq \{0\}$\\
\textbf{(9)}\quad If $I$ is not maximal then $\hat{L}_{I} = \hat{L}^+_{I} \oplus \hat{L}^-_{I}$ where $\hat{L}^{\pm}_{I}$ are spanned by Weyl spinors of $\pm$ chirality and $\dim(\hat{L}^+_{I})=\dim(\hat{L}^-_{I})$.\\
\end{theorem}

Before proving the above results it is important to note that given an isotropic subspace $I\subset V \oplus V^*$ of dimension $k$ then it is always possible to find a vector subspace $V'= \textrm{Span}\{e_1,\ldots,e_n\}$  such that $I$ is spanned by $\{e_1,e_2, \ldots,e_k\}$ and $V \oplus V^*= V' \oplus V'^* $. Hence, the judicious and convenient choice $I= \textrm{Span}\{e_1,\ldots,e_k\}$ represents no loss of generality.
\\
\\
\small{
\emph{Proof of Theorem \ref{Theor-Algebraic}:}\\
\emph{(1)} If $I=\textrm{Span}\{e_1,e_2, \ldots,e_k\}$ then it is not hard to conclude that
$$\hat{L}_I = \textrm{Span}\{\,\hat{1},\, \hat{\theta}^{\alpha'},\, \hat{\theta}^{\alpha'\beta'},\, \cdots,\, \hat{\theta}^{k+1\, k+2... n}\,\}\,,$$
where $\alpha', \beta'\,\in\{k+1,k+2,...,n\}$. Leading us to the following result:
$$\dim(\hat{L}_I)= 1+ (n-k)+ \frac{(n-k)(n-k-1)}{2!}+ \cdots + 1 = \sum_{i=0}^{n-k}\,\dbinom{n-k}{i} = (1+1)^{n-k}. $$
\\
\emph{(2)} Suppose that $\hat{L}_I=\hat{L}_{I'}$ but $I\neq I'$, then there exists $u\in I$ such that $u$ does not belong to $I'$ (or the converse). Now, if $\hat{\varphi}\in \hat{L}_I=\hat{L}_{I'}$ is a non-zero spinor then $2\la u,v'\ra \hat{\varphi}=(uv'+v'u)\cdot \hat{\varphi} = 0$ $\,\forall\,\, v'\in I'$, so that $\la u,v'\ra = 0$ for all $v'\in I'$. This implies that $I'' \equiv I' + \textrm{Span}\{u\}$ is an isotropic subspace. Moreover, since $u\notin I'$ it follows that $\dim(I'')=\dim(I')+1$, from which we conclude that $\dim(\hat{L}_{I''})<\dim(\hat{L}_{I'})$. On the other hand, since  $\hat{L}_I=\hat{L}_{I'}$  then $u\cdot\hat{L}_{I'}=0$, so if $\hat{\varphi}\in \hat{L}_{I'}$ then $v'' \cdot \hat{\varphi} = 0$ for all $v''\in I''$. Thus, $\hat{L}_{I'}\subset\hat{L}_{I''}$ which implies $\dim(\hat{L}_{I''})\geq\dim(\hat{L}_{I'})$, contradicting the former inequality. Therefore, if $\hat{L}_I=\hat{L}_{I'}$ then $I=I'$. The converse is trivial.\\
\\
\emph{(3)} If $\hat{\psi}\in(\hat{L}_{I} + \hat{L}_{I'})$ then $\hat{\psi} = \hat{\varphi} +\hat{\varphi}'$, where $v\cdot \hat{\varphi} = 0 =v'\cdot \hat{\varphi}'$ for all $v\in I$ and $v'\in I'$. Thus, if $u\in I\cap I'$ then $u\cdot\hat{\varphi}= 0 = u\cdot \hat{\varphi}'$, which implies $u\cdot\hat{\psi}=0$. This means that $\hat{\psi}\in\hat{L}_{I\cap I'}$, proving the wanted relation.\\
\\
\emph{(4)} Suppose that $I''=I+I'$ is an isotropic subspace and let $v''\in I''$, then $v'' = v+v'$ where $v\cdot\hat{\varphi}=0=v'\cdot\hat{\varphi}'$ for all $\hat{\varphi}\in\hat{L}_{I}$ and $\hat{\varphi}'\in\hat{L}_{I'}$. Thus, if $\hat{\psi}\in\hat{L}_{I}\cap\hat{L}_{I'}$ then $v''\cdot\hat{\psi}=0$. This means that $\hat{L}_{I}\cap\hat{L}_{I'}\subset\hat{L}_{I''}$. Conversely, if $\hat{\psi}\in\hat{L}_{I''}$ then $(v+v')\cdot\hat{\psi}=0$ for all $v\in I$ and $v'\in I'$. In particular, if we take $v'=0$ then we get $v\cdot\hat{\psi}=0$ for all $v\in I$ $\Rightarrow$ $\hat{\psi}\in\hat{L}_{I}$. In an analogous fashion it is obtained that $\hat{\psi}\in\hat{L}_{I'}$, so that $\hat{\psi}\in \hat{L}_{I}\cap\hat{L}_{I'}$. Therefore, we get  $\hat{L}_{I''}\subset\hat{L}_{I}\cap\hat{L}_{I'}$, establishing the equality $\hat{L}_{I''}=\hat{L}_{I}\cap\hat{L}_{I'}$.\\
\\
\emph{(5)} Suppose that $u\in(V \oplus V^*)$ is such that $u\cdot \hat{\varphi}=0$ for all $\hat{\varphi} \in \hat{L}_{I}$, then if $v\in I$ it follows that $2\la u,v\ra \hat{\varphi} = (uv+vu)\cdot \hat{\varphi} = 0$. Hence, $\la u,v\ra= 0$ for all $v$, so that $I''=I+ \textrm{Span}\{u\}$ is an isotropic subspace of $V \oplus V^*$. By hypothesis it follows that $\hat{L}_{I}\subset\hat{L}_{I''}$. On the other hand, by definition, we have $\dim(I'')\geq\dim(I)$ which implies, by part \emph{(1)} of this theorem, that $\dim(\hat{L}_{I''})\leq\dim(\hat{L}_{I})$. Therefore, we conclude that $\hat{L}_{I}=\hat{L}_{I''}$. Now, using  part \emph{(2)} of this theorem we get $I=I''\equiv I+ \textrm{Span}\{u\}$, thus $u\in I$.\\
\\
\emph{(6)} If $I'\subset I$ then $v'\cdot\hat{\varphi}= 0$ for all $v'\in I'$ and $\hat{\varphi}\in\hat{L}_{I}$, which is tantamount to  $\hat{L}_{I}\subset\hat{L}_{I'}$. Conversely, if $\hat{L}_{I}\subset\hat{L}_{I'}$ then $v'\cdot\hat{\varphi}= 0$ for all $\hat{\varphi}\in\hat{L}_{I}$ and  $v'\in I'$, which thanks to the item \emph{(5)} of this theorem, implies that $v'\in I$, and so $I'\subset I$.\\
\\
\emph{(7)} As explained above, if $\dim(I)=k>0$ then one can assume, without loss of generality, that $I=\textrm{Span}\{e_1,e_2, \ldots,e_k\}$. With this assumption $\hat{L}_I$ is generated by the basis $\{\hat{1},\hat{\theta}^{\alpha'},\hat{\theta}^{\alpha'\beta'},\ldots, \hat{\theta}^{k+1\, k+2... n}\}$, where $\alpha', \beta'\,\in\{k+1,k+2,...,n\}$. The inner product of two elements of this basis, $(\hat{\varphi}_1,\hat{\varphi}_2)$, is always zero, since  $<\hat{\varphi}_1^{\,t}\wedge\hat{\varphi}_2>_n = \pm <\hat{\varphi}_1\wedge\hat{\varphi}_2>_n = 0$. With the last equality stemming from the fact that neither $\hat{\varphi}_1$ nor $\hat{\varphi}_2$ contains the term $\hat{\theta}^{1}$, so that it is impossible to get the spinor $\hat{\theta}^{12...n}$ after the wedge product. Then, using (\ref{InnerSpinors}) we see that $(\hat{\varphi}_1,\hat{\varphi}_2)=0$. The rest of the statement follows from the bilinearity of this inner product.\\
\\
\emph{(8)} If $I''=I+I'$ is isotropic then, from parts \emph{(1)} and \emph{(4)} of this theorem, it follows that $\dim(\hat{L}_{I} \cap\hat{L}_{I'}) = \dim(\hat{L}_{I''}) \geq 1$, so $\hat{L}_{I} \cap\hat{L}_{I'}\neq \{0\}$. Conversely, if $\hat{L}_{I} \cap\hat{L}_{I'}\neq \{0\}$ then there exists $\hat{\varphi}\neq0$ such that $v\cdot\hat{\varphi}=0=v'\cdot\hat{\varphi}$ for all $v\in I$ and $v'\in I'$. Thus, $2\la v,v'\ra\hat{\varphi} = (vv'+v'v)\cdot\hat{\varphi}=0$, so that $\la v,v'\ra=0$ for all $v\in I$ and $v'\in I'$, implying that $I''=I+I'$ is isotropic.\\
\\
\emph{(9)} This result is a simple consequence of the fact that if $I$ is non-maximal, $k\neq n$, then it can always be chosen to be $I=\textrm{Span}\{e_1,e_2, \ldots,e_k\}$, so that $\hat{L}_I = \textrm{Span}\{\hat{1},\hat{\theta}^{\alpha'},\hat{\theta}^{\alpha'\beta'},\ldots, \hat{\theta}^{k+1\, k+2... n}\}$, where $\alpha', \beta'\,\in\{k+1,k+2,...,n\}$. In this basis for $\hat{L}_I$ there are $2^{n-k-1}$ spinors of even grade and $2^{n-k-1}$ spinors of odd grade. \qed
}\normalsize
\\

A corollary that easily follows from items \emph{(1)}, \emph{(4)} and \emph{(8)} is that if two isotropic subspaces $I$ and $I'$ are such that the dimension of $\hat{L}_{I}\cap\hat{L}_{I'}$ is not even then $\dim(\hat{L}_{I}\cap\hat{L}_{I'}) =1$ and $(I+I')$ is a maximally isotropic subspace. Furthermore, part  \emph{(3)} of Theorem 1 states that $(\hat{L}_{I} + \hat{L}_{I'})$ is contained in the space  $\hat{L}_{I\cap I'}$. Then, a natural question to be posed is whether this result can be made stronger. For instance, is it true that both spaces are equal? The answer is no, generally there are elements of $\hat{L}_{I\cap I'}$ that are not contained in $(\hat{L}_{I} + \hat{L}_{I'})$. More precisely, it can be proved that if $(\hat{L}_{I} + \hat{L}_{I'})= \hat{L}_{I\cap I'}$ then one of the special cases must happen: (i) $I\subset I'$, (ii) $I'\subset I$ or (iii) $\dim(I)=\dim(I')=\dim(I\cap I')+1$ with $(I+I')$ non-isotropic.

It is not hard to see that the pseudo-scalar $\mathcal{I}\in Cl(V\oplus V^*)$ anti-commutes with every vector $u \in V\oplus V^*$, $\mathcal{I}u= -u\mathcal{I}$. Thus, if $I\subset V\oplus V^*$ is some isotropic subspace and $\hat{\varphi}\in \hat{L}_I$ then $v\cdot(\mathcal{I}\cdot\hat{\varphi})= - \mathcal{I}\cdot(v\cdot\hat{\varphi})= 0$ for all $v\in I$. Therefore, the spinor $(\mathcal{I}\cdot\hat{\varphi})$ must belong to $\hat{L}_I$. In particular, if $\hat{\varphi}\in \hat{L}_I$ then the Weyl spinors of positive and negative chirality $\hat{\varphi}^\pm= \frac{1}{2}(\hat{\varphi}\pm \mathcal{I}\cdot\hat{\varphi})$ belong to $\hat{L}_I$ as well. This, in turn, implies that $\hat{L}_I$ can be split as the direct sum of a subspace of positive chirality and a subspace of negative chirality, which agrees with item \emph{(9)} of Theorem \ref{Theor-Algebraic}. Particularly, such reasoning entails the well-known result that every pure spinor must be a Weyl spinor. It is also worth noting that part \emph{(7)} of this theorem guarantees that every pure spinor $\hat{\psi}$ must be orthogonal to itself, $(\hat{\psi}, \hat{\psi})= 0$.

Another interesting problem that can be posed is the following: Given a proper subspace of the spinor space, $\Upsilon\subset S$, then how can one know whether $\Upsilon$ it is related to some isotropic subspace $I\subset V \oplus V^*$ ? More precisely, given $\Upsilon$ how to know if there exists some isotropic subspace $I$ such that $\Upsilon=\hat{L}_{I}$? The above theorem provides several clues for the answer: If one of the conditions \textit{(i) $\dim(\Upsilon)\neq 2^s$ for some $s\in\{0,1,...,n-1\}$, (ii)  $(\hat{\varphi}_1,\hat{\varphi}_2)\neq0$ for some  $\hat{\varphi}_1,\hat{\varphi}_2\,\in\Upsilon$, (iii) $\dim(\Upsilon)=1$ and $\hat{\varphi}\in\Upsilon$ is not a Weyl spinor} or \textit{(iv) $\dim(\Upsilon)\neq1$ and $\Upsilon$ does not admit a basis such that half of its elements are Weyl spinors of positive chirality while the other half are Weyl spinors of negative chirality} happen then the spinorial subspace $\Upsilon$ is not related to any isotropic subspace. But the converse is not true. For instance, if $\Upsilon$ is the space spanned by some $\hat{\varphi}^+\in S^+$ then in general $\Upsilon$ is not related to an isotropic subspace, since it is well-known that although every pure spinor is a Weyl spinor not all Weyl spinors are pure\footnote{If $n\leq3$ then all Weyl spinors are pure. But if $n>3$ there are more Weyl spinors than pure spinors.}. Then, it is convenient to make the following definition:\\
\\
\textbf{Definition:} \emph{A proper spinorial subspace $\Upsilon\subset S$ is called pure if there exists some isotropic subspace $I\subset V \oplus V^*$, such that $\Upsilon=\hat{L}_{I}$. In the latter case we shall say that $\Upsilon$ is the pure subspace associated to $I$.}
\\
\\
Since the constraint that a spinor must obey in order for it to be pure is quadratic in the spinor \cite{Cartan,Spinor-thesis}, then, probably, the algebraic conditions that a proper spinorial subspace $\Upsilon\subset S$ might obey in order for it to be a pure subspace might consist of quadratic equations in the spinors of $\Upsilon$.

As a last comment in this section, note that the space $\hat{L}_{I}$ can be defined in a different way from the one presented in Eq. (\ref{Spinor subspace}). When $I$ is maximal then $\hat{L}_{I}$ is the pure spinor line related to $I$. When $I$ is not maximal there are several ways to complete $I$ in order to create maximally isotropic subspaces. Now, let $\{I_a\}$ be the set of all isotropic subspaces of dimension $n$ such that $I\subset I_a$. Associated to each maximally isotropic subspace $I_a$ there is a pure spinor line spanned by $\hat{\varphi}_a$. Then, the space $\hat{L}_{I}$ is just the one spanned by the pure spinors $\{\hat{\varphi}_a\}$.
\\
\\
\textbf{Remark:}\emph{ It is not difficult to note that if $\hat{\varphi}_1, \hat{\varphi}_2\in\hat{L}_I$ and $(e+\theta)\in I$ then $(e+2\theta)\cdot(\hat{\varphi}_1\wedge \hat{\varphi}_2)=0$. In particular, this implies that if $I'=\textrm{Span}\{e_1,e_2, \ldots,e_k\}$ then $(\hat{\varphi}'_1\wedge \hat{\varphi}'_2)\in\hat{L}_{I'}$ for all $\hat{\varphi}'_1, \hat{\varphi}'_2\in\hat{L}_{I'}$. But this is not true in the case of $I'$ being an arbitrary isotropic subspace. So, probably there exists some bilinear operation $\diamond:S\times S\rightarrow S$ such that given a general isotropic subspace $I$ then $(\hat{\varphi}_1\diamond \hat{\varphi}_2)$ belongs to $\hat{L}_{I}$ for all $\hat{\varphi}_1, \hat{\varphi}_2 \in \hat{L}_{I}$ and such that $\diamond$ reduces to $\wedge$ when $I$ is equal to $I'$. }

\section{Manifolds and Spinorial Connections}\label{section-connection}

From now on we are going to work on a $2n$-dimensional manifold $(M,g)$ endowed with a non-degenerate metric $g$. Furthermore, the tangent bundle $TM$ is assumed to be endowed with a torsion-free and metric-compatible derivative, the Levi-Civita connection. Since in the present article we are concerned only with local results it follows that we are allowed to identify the complexified tangent spaces $\mathbb{C}\otimes T_pM$, at any point $p\in M$, with a vector space of the form $V\oplus V^*$, so that the results of the previous sections can be used. More precisely, if $\{e_i,\theta^j\}$ is a frame in a patch of $M$, with $i,j\in \{1,2,\ldots,n\}$, such that
$$ g(e_i,e_j)\,=\, 0\,, \quad\;  g(e_i,\theta^j)\,=\, \frac{1}{2}\,\delta^j_i\,, \quad\;   g(\theta^i,\theta^j)\,=\, 0\,,   $$
then we shall say that $V$ is the vector bundle spanned by $\{e_i\}$ while $V^*$ is the dual vector bundle spanned by $\{\theta^i\}$, with $\theta^j(e_i)\equiv \delta^j_i$. So that an equation analogous to (\ref{Inner Product}) holds:
$$   g( e+\theta, e'+ \theta') \,\equiv\, \frac{1}{2}\left[\,\theta(e') \,+\, \theta'(e)\,\right]\;\;\;\;\;\; e,e'\in \Gamma(V) \;\textrm{ and }\;  \theta, \theta'\in \Gamma(V^*) \,.  $$
Where $\Gamma(V)$ is the space of sections of the bundle $V$. So, in this approach a fiber of the spinorial bundle $S$ will be identified with $\wedge V^*$. It is worth remarking that the identifications $\mathbb{C}\otimes T_pM\sim V|_p\oplus V^*|_p$ and $S_p\sim \wedge V^*|_p$ are not canonical at all, they depend on an arbitrary choice of frame. Generally, the frame $\{e_i,\theta^j\}$ cannot be globally defined (non-parallelizable manifolds), so that these objects are well-defined just in a local neighborhood of the arbitrary point $p\in M$. Therefore, the space $V$ might, actually, be understood as a local trivialization of some bundle, which is sufficient for our purposes.

Before proceeding, let us establish some notation. It is useful to denote the frame $\{e_i,\theta^j\}$ by $\{e_a\}$, where the indices $a,b,\ldots$ run from $1$ to $2n$ and $e_{i+n}\equiv \theta^i$. The components of the metric in this frame are denoted by the symmetric matrix $g_{ab}\equiv g(e_a,e_b)$, while $g^{ab}$ stands for its inverse, $g^{ab}g_{bc}= \delta^a_c$. Finally, let $\omega_{ab}^{\ph{ab}c}$ be the connection coefficients in this frame:
$$  \nabla_a\,e_b \,\equiv\, \nabla_{e_a}\,e_b \,=\, \omega_{ab}^{\ph{ab}c}\,e_c \,. $$
It is customary to raise or lower the indices by means of $g^{ab}$ and $g_{ab}$ respectively. For instance, $\omega_{abd}= \omega_{ab}^{\ph{ab}c}g_{cd}$ and $\omega_{a}^{\ph{a}dc}= \omega_{ab}^{\ph{ab}c}g^{bd}$. Since the components of the matrix $g_{ab}$ are constant and the connection is metric compatible, it follows that $\omega_{abc}= -\omega_{acb}$.

Now, let us introduce a connection $\hat{\nabla}_a$ on the spinorial bundle $S$. In order for this connection to be an extension of the Levi-Civita connection one might impose for it to satisfy the Leibniz rule with respect to the Clifford action:
\begin{equation}\label{Leibniz-Cliff}
 \hat{\nabla}_a\,(v\cdot \hat{\varphi}) \,=\, (\nabla_a\,v)\cdot \hat{\varphi} \,+\, v\cdot\,\hat{\nabla}_a\,\hat{\varphi}\quad\;\; \forall\;\; v\in\Gamma(TM) \;,\,\,\hat{\varphi}\in \Gamma(S)\,.
 \end{equation}
Then, let $\{\hat{\psi}_\alpha\}$ be a local frame of the spinorial bundle $S$, with $\alpha,\beta,\ldots$ running from $1$ to $2^n$, such that the Clifford action of the frame $\{e_a\}$ on the spinors $\hat{\psi}_\alpha$ is constant in the relevant patch of $M$. This means that the following relation holds:
\begin{equation}\label{e-constant}
 e_a\cdot \hat{\psi}_\alpha \,\equiv\, (e_a)^\beta_{\ph{\beta}\alpha} \,\hat{\psi}_\beta\,, \quad \textrm{with the matrices }  \;(e_a)^\beta_{\ph{\beta}\alpha} \; \textrm{ being constant}.
\end{equation}
In physics the matrices $(e_a)^\beta_{\ph{\beta}\alpha}$ are the so-called Dirac matrices. For instance, we can assume that the frame $\{\hat{\psi}_\alpha\}$ is given by $\{\hat{1}, \hat{\theta^i}, \hat{\theta}^{ij}, \ldots, \hat{\theta}^{1\ldots n}\}$. Then, without loss of generality, define the action of the connection on this frame of spinors to be:
\begin{equation}\label{D-SpinorBasis}
\hat{\nabla}_a\, \hat{\psi}_\alpha \,\equiv\, \Omega_a\cdot\hat{\psi}_\alpha \,=\, (\Omega_a)^\beta_{\ph{\beta}\alpha} \,\hat{\psi}_\beta \,.
\end{equation}
Now, we might look for the general expression of $\Omega_a$ that satisfies Eq. (\ref{Leibniz-Cliff}). Computing $\hat{\nabla}_a(e_b\cdot\hat{\psi}_\alpha)$ in one hand by means of (\ref{e-constant}) and on the other hand by means of (\ref{Leibniz-Cliff}) and then equating both calculations lead us to the following result:
$$  (e_b)^\beta_{\ph{\beta}\alpha}\, (\Omega_a)^\sigma_{\ph{\sigma}\beta} \,\hat{\psi}_\sigma \,\,=\,\,
  \omega_{ab}^{\ph{ab}c}\,(e_c)^\sigma_{\ph{\sigma}\alpha}\,\hat{\psi}_\sigma \,+\,
   (\Omega_a)^\beta_{\ph{\beta}\alpha} \, (e_b)^\sigma_{\ph{\sigma}\beta}\,\hat{\psi}_\sigma \,. $$
In terms of Clifford algebra this equation is tantamount to:
\begin{equation}\label{Omega-relation}
  \omega_{ab}^{\ph{ab}c}\,e_c \,=\,   \Omega_a\,e_b \,-\, e_b \,\Omega_a\,.
\end{equation}
Such relation makes clear that $\Omega_a$ is defined up to the sum of a term on the centre of the Clifford algebra. But when the dimension is even, as assumed throughout this paper, the centre of the Clifford algebra is spanned by the identity operator. Then, the general solution to Eq. (\ref{Omega-relation}) is:
\begin{equation}\label{Omega-Def}
  \Omega_a \,=\, -\,\frac{1}{4} \omega_{a}^{\ph{a}bc}\, (e_b\wedge e_c)\,+\, A_a \,=\, -\,\frac{1}{4} \omega_{a}^{\ph{a}bc}\, e_b\,e_c \,+\, A_a \,.
\end{equation}
Where $A_a$ are the components of an arbitrary $1$-form $A\in \Gamma(TM^*)$. Thus, if $v= v^b e_b$ and $\hat{\varphi}= \varphi^\beta \hat{\psi}_\beta$  are general vector and spinor fields respectively, then Eqs. (\ref{Omega-relation}) and (\ref{D-SpinorBasis}) imply that:
\begin{equation}\label{Connection-Spinor}
\nabla_{a}\, v \,=\, (\partial_a v^b)e_b + v^b(\Omega_a e_b- e_b\Omega_a) \quad;\quad \hat{\nabla}_a\,\hat{\varphi} \,=\,
(\partial_a\varphi^\beta)\hat{\psi}_\beta +  \varphi^\beta\, \Omega_a\cdot   \hat{\psi}_\beta \,\,.
\end{equation}
A nice review about connections on the spinorial bundle is available in \cite{Trautman-connection}, see also \cite{Benn-Book}.

Note that since the spinorial space is defined to be the space where an irreducible and minimal representation of the Clifford algebra acts, it follows that  spinors are defined up to a conformal scale. For instance, if the frame $\{\hat{\psi}_\alpha\}$ obeys to Eq. (\ref{e-constant}) then the frame formed by the spinors $\hat{\psi}'_\alpha\equiv e^{-\lambda}\hat{\psi}_\alpha$,  for some function $\lambda$, also obeys to the same equation. Likewise, it is worth recalling that the inner product (\ref{InnerSpinors}) was also defined up to a global multiplicative factor, so that one can easily redefine this inner product, $(\,,)\rightarrow \, \prec\,,\succ$, in such a way that  $(\hat{\psi}_\alpha, \hat{\psi}_\beta)= \prec\hat{\psi}'_\alpha, \hat{\psi}'_\beta\succ$. Thus one can see the transformation $\hat{\varphi}\mapsto \hat{\varphi}'= e^{-\lambda}\hat{\varphi}$ as an intrinsic symmetry of the spinorial formalism. However, in order for this symmetry to be compatible with the spinorial connection we must assume that the connection introduced in (\ref{Omega-Def}) and (\ref{Connection-Spinor}) change in a way that $A_a\mapsto A'_a\,=\, A_a + \partial_a\lambda$. More precisely, we shall define the connection $\hat{\nabla}'_a$ to be such that
$$ \hat{\nabla}'_a\, \hat{\psi}_\alpha \,=\, \Omega'_a\cdot\hat{\psi}_\alpha \,=\, (\Omega_a+ \partial_a\lambda)\cdot\hat{\psi}_\alpha \,=\,  \hat{\nabla}_a\, \hat{\psi}_\alpha + (\partial_a\lambda)\, \hat{\psi}_\alpha\,. $$
With this definition one arrives at the following relation:
$$  \hat{\varphi}' \,=\, e^{-\lambda}\,\hat{\varphi} \; \textrm{ and } A'_a\,=\, A_a + \partial_a\lambda \; \;\;\Rightarrow\quad \hat{\nabla}'_a\,\hat{\varphi}' \,=\, e^{-\lambda} \hat{\nabla}_a\,\hat{\varphi} \,=\, \lef\hat{\nabla}_a\,\hat{\varphi} \rig'\,.
 $$
Physically, the above equation says that $\hat{\nabla}_a$ is the covariant derivative associated with the scaling symmetry of the spinors.

%$$ \hat{\varphi}\,\mapsto\, \hat{\varphi}' \,=\, e^{-\lambda}\,\hat{\varphi} \quad\quad\Rightarrow\quad \quad
%A_a\,\mapsto\, A'_a\,=\, A_a + \partial_a\lambda \,. $$

The curvature of the spinorial bundle is defined by the following action:
$$ \hat{\mathfrak{R}}_{ab}\,\hat{\psi} \,=\, \left( \hat{\nabla}_a\hat{\nabla}_b - \hat{\nabla}_b\hat{\nabla}_a - \hat{\nabla}_{[e_a,e_b]} \right) \,\hat{\psi} \,.$$
Where $[e_a,e_b]$ means the Lie bracket of the vector fields $e_a$ and $e_b$. Then, using (\ref{D-SpinorBasis}) and (\ref{Omega-Def}), one can explicitly prove that this curvature is given by:
\begin{equation}\label{Curvature-Action}
 \hat{\mathfrak{R}}_{ab}\,\hat{\psi} \,=\, -\frac{1}{4}\, R_{ab}^{\phantom{ab}cd}\,(e_c e_d)\cdot\hat{\psi} \,+\, F_{ab}\,\hat{\psi}\;\,\,.
\end{equation}
Where, in the above equation, $R_{ab}^{\phantom{ab}cd}= g^{cf}g^{dh}R_{abfh}$ is the Riemann tensor of the Riemannian manifold $(M,g)$, while $F_{ab}$ are the components of the exterior derivative of the $1$-form $A$, $F_{ab}= (dA)_{ab}$. In particular, note that the spinorial curvature is invariant under the scaling transformation, $A_a\mapsto A'_a\,=\, A_a + \partial_a\lambda$.

Two other operators that are worth mentioning are the Dirac and the twistor operators, defined respectively by:
$$ \hat{\bl{D}}\,=\, g^{ab}e_a\cdot\hat{\nabla}_b \quad \;;\;\quad \hat{\bl{T}}_a \,=\, \hat{\nabla}_a \,-\, \frac{1}{2n}\,e_a \cdot\hat{\bl{D}} \,.$$
The twistor operator is characterized by the property $g^{ab}e_a\cdot\hat{\bl{T}}_b= 0$. We shall say that a spinor $\hat{\psi}$ is a twistor if it is annihilated by the action of the twistor operator, $\hat{\bl{T}}_a \hat{\psi}= 0$ $\forall\, a$. As an aside, note that the square of the Dirac operator is given by:
$$ \hat{\bl{D}}^2\, \hat{\psi} \,=\, \hat{\Box} \,\hat{\psi} \,+\, \frac{1}{2}F^{ab}\,(e_ae_b)\cdot\, \hat{\psi} \,+\, \frac{1}{4}\,R\, \hat{\psi} \,.$$
Where $\hat{\Box}$ is the ``spinorial Laplacian'', $\hat{\Box}= \hat{\nabla}^a\hat{\nabla}_a- \hat{\nabla}_{\nabla^ae_a}$, and $R=R_{ab}^{\phantom{ab}ab}$ is the Ricci scalar.

So far, in order to define a connection on the spinorial bundle we have just imposed that this connection obeys to the Leibniz rule with respect to the Clifford action, see (\ref{Leibniz-Cliff}). As a consequence, we have found that such connection is unique up to an additive $1$-form $A$. It turns out that such freedom can be fixed once we require the spinorial connection to obey the Leibniz rule with respect to the inner product on the spinorial bundle,  in addition to Eq. (\ref{Leibniz-Cliff}). More precisely, if one requires that
\begin{equation}\label{Leibniz-InnerProd}
 \nabla_a\, (\hat{\varphi},\hat{\psi}) \,=\, (\hat{\nabla}_a\hat{\varphi},\hat{\psi}) \,+\, (\hat{\varphi},\hat{\nabla}_a\hat{\psi})
\end{equation}
then the choice of $A_a$ in Eq. (\ref{Omega-Def}) is unique. For example, if we choose the spinorial frame $\{\hat{\psi}_\alpha\}$ to be $\{\hat{1}, \hat{\theta^i}, \hat{\theta}^{ij}, \ldots, \hat{\theta}^{1\ldots n}\}$ then Eq. (\ref{Leibniz-InnerProd}) holds if, and only if, we set $A_a= 0$.
For instance, this was the choice of connection made by R. Penrose when he introduced the spinorial formalism in 4-dimensional general relativity \cite{Penrose-Book}, in index notation Eq. (\ref{Leibniz-InnerProd}) means that the symplectic form $\epsilon_{AB}$ is covariantly constant. For sake of generality, in the following sections it will not be assumed that Eq. (\ref{Leibniz-InnerProd}) holds, so that the freedom in the choice of $A_a$ can be exploited.

%%%%%%%%%%%%%%%
%%Now if $\hat{\varphi}= f^\alpha\hat{\psi}_\alpha$ is a general spinor field the we can define the following connection on $S$:
%%$$ \hat{\nabla}_a\,\hat{\varphi}\,=\, \lef\partial_a\,f^\alpha\rig\hat{\psi}_\alpha  \,-\,\frac{1}{4} f^\alpha \, \omega_{a}^{\ph{a}bc}\,
%%(e_b\wedge e_c)\cdot \hat{\psi}_\alpha$$
%%%%%%%%%%%%%%%%%%%%

%%%%%%%%%%%%%%%%%%%%%%%%%%%%%%%%%%%%%%%%%%%%%%%%%%%%%%%%%%%%%%%%%%%%%%%%%%%%%%%%%%%%%%%%%%%%%%%%%%%%%
%%%%%%%%%%%%%%%%%%%%%%%%%%%%%%%%%%%%%%%%%%%%%%%%%%%%%%%%%%%%%%%%%%%%%%%%%%%%%%%%%%%%%%%%%%%%%%%%%%%%%
%%%%%%%%%%%%%%%%%%%%%%%%%%%%%%%%%%%%%%%%%%%%%%%%%%%%%%%%%%%%%%%%%%%%%%%%%%%%%%%%%%%%%%%%%%%%%%%%%%%%%
%%%%%%%%%%%%%%%%%%%%%%%%%%%%%%%%%%%%%%%%%%%%%%%%%%%%%%%%%%%%%%%%%%%%%%%%%%%%%%%%%%%%%%%%%%%%%%%%%%%%%
%%%%%%%%%%%%%%%%%%%%%%%%%%%%%%%%%%%%%%%%%%%%%%%%%%%%%%%%%%%%%%%%%%%%%%%%%%%%%%%%%%%%%%%%%%%%%%%%%%%%%
%%%%%%%%%%%%%%%%%%%%%%%%%%%%%%%%%%%%%%%%%%%%%%%%%%%%%%%%%%%%%%%%%%%%%%%%%%%%%%%%%%%%%%%%%%%%%%%%%%%%%
%%%%%%%%%%%%%%%%%%%%%%%%%%%%%%%%%%%%%%%%%%%%%%%%%%%%%%%%%%%%%%%%%%%%%%%%%%%%%%%%%%%%%%%%%%%%%%%%%%%%%
%%%%%%%%%%%%%%%%%%%%%%%%%%%%%%%%%%%%%%%%%%%%%%%%%%%%%%%%%%%%%%%%%%%%%%%%%%%%%%%%%%%%%%%%%%%%%%%%%%%%%

\section{Pure Subspaces and Integrability}\label{section-Integrability}

Once we have introduced a connection on the spinorial bundle, one can look for theorems on the integrability of isotropic distributions and its relation with the pure subspaces. For instance, let $I$ be some maximally isotropic distribution of vector fields over $(M,g)$ and $\hat{\psi}$ the associated pure spinor, meaning that $\hat{\psi}$ is annihilated by $I$. Then, it is well-known that the distribution $I$ is integrable if, and only if, $\hat{\nabla}_X\hat{\psi}\propto \hat{\psi}$ for all $X\in I$. Using the formalism introduced in Section \ref{section-AlgebRes}, it turns out that this result can be seen as a particular case of a more general result stated here in the form of the following theorem.
\begin{theorem}\label{Theo. Integrab}
An isotropic distribution of vector fields $I$ is integrable if, and only if, $X\cdot\hat{\nabla}_Y\hat{\varphi}= Y\cdot\hat{\nabla}_X\hat{\varphi}$ for all $X,Y\in I$ and $\hat{\varphi} \in \hat{L}_{I}$.
\end{theorem}
\vspace{0.3cm}
\small{
\emph{Proof of Theorem \ref{Theo. Integrab}:}\\
Let $I= I_k= \textrm{Span}\{e_1,e_2,...,e_k \}$, with $k$ fixed, and $\hat{\varphi} \in \hat{L}_{I_k}$. Then, $e_\beta\cdot\hat{\varphi}=0$, where $\beta\in\{1,2,...,k\}$. Hence,
\begin{align}
   \nonumber 0=&\hat{\nabla}_\alpha(e_\beta\cdot\hat{\varphi})= (\nabla_\alpha e_\beta)\cdot\hat{\varphi} + e_\beta\cdot(\hat{\nabla}_\alpha\hat{\varphi}) \;\;\Rightarrow \\
   &[e_\alpha,e_\beta]\cdot\hat{\varphi} = e_\alpha\cdot(\hat{\nabla}_\beta\hat{\varphi}) - e_\beta\cdot(\hat{\nabla}_\alpha\hat{\varphi})\,.\label{Theorem2}
\end{align}
Now, if $I$ is integrable then $[e_\alpha,e_\beta] = f_{\alpha\beta}^{\gamma}e_\gamma$, where $\alpha,\beta,\gamma\in\{1,2,...,k\}$, so that $[e_\alpha,e_\beta]\cdot\hat{\varphi} = 0$. Conversely, by item \emph{(5)} of Theorem \ref{Theor-Algebraic}, if $[e_\alpha,e_\beta]\cdot\hat{\varphi} = 0$ for all $\hat{\varphi} \in \hat{L}_{I_k}$ then $[e_\alpha,e_\beta] \in I$, so $I$ is integrable. Therefore, (\ref{Theorem2}) implies that $I_k$ is integrable $\Leftrightarrow$ $e_\alpha\cdot(\nabla_\beta\hat{\varphi}) = e_\beta\cdot(\nabla_\alpha\hat{\varphi})$, for all $\alpha,\beta\in\{1,2,...,k\}$ and for all $\hat{\varphi} \in \hat{L}_{I_k}$, proving Theorem 2.\qed
}\normalsize\\

Now, since $I_k$ is an isotropic distribution it follows that $g(e_\alpha,e_\beta)=0$, which in terms of Clifford algebra means that $e_\alpha e_\beta = -e_\beta e_\alpha$. Then, supposing that $I_k$ is integrable and using the last theorem we have,
\begin{align*}
    (e_\alpha e_\beta)\cdot \hat{\nabla}_\gamma\hat{\varphi} \,&=\, -(e_\beta e_\alpha)\cdot \hat{\nabla}_\gamma\hat{\varphi} \,=\, - (e_\beta e_\gamma)\cdot \hat{\nabla}_\alpha\hat{\varphi} \,=\,  (e_\gamma e_\beta)\cdot \hat{\nabla}_\alpha\hat{\varphi} \\
    \,&=\, (e_\gamma e_\alpha)\cdot \hat{\nabla}_\beta\hat{\varphi} \,=\, - (e_\alpha e_\gamma)\cdot \hat{\nabla}_\beta\hat{\varphi} \,=\, - (e_\alpha e_\beta)\cdot \hat{\nabla}_\gamma\hat{\varphi}\,,
\end{align*}
with $\hat{\varphi} \in \hat{L}_{I_k}$ and $\alpha,\beta,\gamma\in\{1,2,...,k\}$. This means that
$v\cdot(e_\beta\cdot \hat{\nabla}_\gamma\hat{\varphi})=0$ for all $v\in I_k$. Hence, by definition, we conclude that $(e_\beta\cdot \hat{\nabla}_\gamma\hat{\varphi})$ belongs to $\hat{L}_{I_k}$, proving the following corollary.
\begin{corollary}
If $I$ is an integrable isotropic distribution then $[X\cdot(\hat{\nabla}_Y\hat{\varphi})]$ belongs to $\hat{L}_{I}$ for all $X,Y\in I$ and  $\hat{\varphi} \in \hat{L}_{I}$.
\end{corollary}

In the case of $I_k$ being maximal we have $I_k=I_n=\textrm{Span}\{e_1,e_2,...,e_n \}$, so that $\hat{L}_{I_k}$ is generated by the pure spinor $\hat{1}$. So, the above corollary implies that if $I_n$ is integrable then $e_i\cdot(\hat{\nabla}_j\hat{1})\propto\hat{1}$ for all $i,j\in\{1,2,...,n\}$, which lead us to the relation $\hat{\nabla}_j\hat{1}= \lambda_j \hat{1}+\kappa_{ji}\hat{\theta}^i$. But the covariant derivative of a chiral spinor cannot change its chirality, so that we must have $\kappa_{ji}=0$. Hence, if $I$ is a maximally isotropic distribution that is integrable then the parallel transport of its  pure spinor in a direction tangent to $I$ does not change the direction of the spinor. By Theorem \ref{Theo. Integrab}, it is clear that the converse of this result is also valid. This is a well-known result that in the language of the pure subspaces assumes the following form.
\begin{corollary}\label{Coro-MaxIsoInteg}
A maximally isotropic distribution $I$ is integrable if, and only if, $\hat{\nabla}_X\,\hat{\varphi}\propto\hat{\varphi}$ for all $X\in I$ and $\hat{\varphi}\in \hat{L}_{I}$.
\end{corollary}

Now, let $I_1=\textrm{Span}\{e_1\}$ be some one-dimensional distribution generated by the null vector field $e_1$, and $\hat{L}_{I_1}$ its associated pure subspace. So, if $\hat{\varphi}\in\hat{L}_{I_1}$ then  $(\nabla_1e_1) \cdot \hat{\varphi} = - e_1\cdot (\hat{\nabla}_1\hat{\varphi})$. Such relation, along with item \emph{(5)} of Theorem \ref{Theor-Algebraic}, implies that the vector field $e_1$ is geodesic if, and only if, $(\hat{\nabla}_1\hat{\varphi})$ belong to $\hat{L}_{I_1}$ for all $\hat{\varphi}\in\hat{L}_{I_1}$. This simple result is just a particular case of a broader theorem concerning the totally geodesic character of isotropic foliations. Before stating and proving such theorem, let us recall some properties of totally geodesic submanifolds.

If $M'$ is a submanifold of $(M,g)$ then a point $p\in M'\subset M$ is said to be geodesic when every geodesic of $M$ that is tangent to $M'$ at the point $p$ remains in $M'$ forever. The submanifold $M'$ is said to be totally geodesic if all its points are geodesic. As an example, note that a geodesic curve in $M$ is always a one-dimensional totally geodesic submanifold of $M$. In order to proceed it is important to explicitly show which restrictions are imposed to the connection coefficients by the existence of a totally geodesic submanifold. Let $\{x^1,x^2,...,x^{2n}\}$ be local coordinates for $M$ in the neighborhood of $p\in M'\subset M$ such that $\{x^1,x^2,...,x^{k}\}$ are local coordinates for $M'$ in this neighborhood. Thus, $\{\partial_1,\partial_2,...,\partial_k\}$ spans the tangent spaces of $M'$ near $p$. Now, let $x^\mu(t)$ be a geodesic of $M$ such that $x(0)= p$ and $\frac{dx}{dt}(0)=\partial_1|_p$. Then, using the geodesic equation it is trivial to see that near $p$
\begin{align}
   \nonumber x^\mu(t) =& x^\mu(0) + t\,\delta_1^{\phantom{1}\mu} - \frac{1}{2}t^2\,\Gamma_{11}^{\phantom{11}\mu} + O(t^3)\;\;\Rightarrow
    \\     &\frac{dx^\mu}{dt}(t)= \delta_1^{\phantom{1}\mu} - t\,\Gamma_{11}^{\phantom{11}\mu} + O(t^2)\label{X(t)}\,.
\end{align}
Where $\Gamma_{\mu\nu}^{\phantom{\mu\nu}\rho}$ is the Christoffel symbol of the metric $g$ evaluated at $p$. Now, if $p$ is a geodesic point then $x^\mu(t)$ must be a point in $M'$ and $\frac{dx^\mu}{dt}(t)$ must be tangent to $M'$ for all $t$. In this case Eq. (\ref{X(t)}) implies that $\Gamma_{11}^{\phantom{11}\mu}|_p= 0$ if $\mu >k$. In general, if $p$ is a geodesic point then $\Gamma_{\alpha\alpha}^{\phantom{11}\mu}|_p =0$ (no sum in $\alpha$) for all $\alpha\in \{1,2,...,k\}$ and $\mu\in \{k+1,k+2,...,2n\}$. It is not hard to note that this is also a sufficient condition. By means of this result along with the Frobenius theorem, we are led to the following conclusion: \emph{The leaves of an integrable distribution $I'$ are totally geodesic submanifolds if, and only if, $\nabla_XY$ is tangent to $I'$ for all $X,Y$ tangent to $I'$}.

There are some other equivalent ways to characterize an integrable distribution of totally geodesic leaves. Let $\{E_1,E_2,...,E_{2n}\}$ be vector fields that form a frame in $(M,g)$ such that $\{E_1,E_2,...,E_{k}\}$ span the leaves of an integrable distribution. Then, these leaves are totally geodesic submanifolds if, and only if,
$$\nabla_X(E_1\wedge E_2\wedge \cdots\wedge E_k)\propto (E_1\wedge E_2\wedge \cdots\wedge E_k) \quad\; \forall\;X\,\in \textrm{Span}\{E_1,E_2,...,E_{k}\} \,. $$
Analogously, if $\{E^1,E^2,...,E^{2n}\}$ is the dual frame of 1-forms, $E^a(E_b)=\delta_b^{\phantom{b}a}$, then the leaves of the integrable distribution are totally geodesic if, and only if,
$$ \nabla_X(E^{k+1}\wedge \cdots\wedge E^{2n})\propto (E^{k+1}\wedge \cdots\wedge E^{2n}) \quad\; \forall\;X\,\in \textrm{Span}\{E_1,E_2,...,E_{k}\} $$
Finally, it is worth remarking that an embedded Euclidean manifold $M'$ is totally geodesic if, and only if, its second fundamental form vanishes. Now, we are ready to prove the following theorem:
\begin{theorem}\label{Theo. - Tot. Geod.}
An isotropic distribution of vector fields $I$ is integrable and its leaves are totally geodesic submanifolds if, and only if,  $(\nabla_X\hat{\varphi})$ belongs to $\hat{L}_{I}$ for all $X\in I$ and $\hat{\varphi} \in \hat{L}_{I}$.
\end{theorem}
\vspace{0.3cm}
\small{
\emph{Proof of Theorem \ref{Theo. - Tot. Geod.}:}\\
If $I$ is integrable and generates totally geodesic submanifolds then, by what was seen above, it follows that $(\nabla_XY)\in I$ for all $X,Y\in I$, so that
$$ 0=(\nabla_XY)\cdot\hat{\varphi}=-Y\cdot(\hat{\nabla}_X\hat{\varphi}) \;\;\;\forall\; X,Y\in I\;\,\textrm{and}\;\,\hat{\varphi} \in \hat{L}_{I}.
$$
Since $Y$ is any vector field of $I$ then, by definition, $(\hat{\nabla}_X\hat{\varphi})$ belongs to $\hat{L}_{I}$. Conversely, if $(\hat{\nabla}_X\hat{\varphi}) \in \hat{L}_{I}$ for all $X\in I$ and  $\hat{\varphi} \in \hat{L}_{I}$ then it follows that $Y\cdot(\nabla_X\hat{\varphi})=0=X\cdot(\hat{\nabla}_Y\hat{\varphi})$. So, by Theorem 2 this implies that the distribution $I$ is integrable. Moreover, since $0=Y\cdot(\nabla_X\hat{\varphi})= - (\nabla_XY)\cdot\hat{\varphi}$ for all $\hat{\varphi} \in \hat{L}_{I}$ then by item \emph{(5)} of Theorem \ref{Theor-Algebraic} we find that $(\nabla_XY)$ must belong to $I$, implying that the leaves of this isotropic distribution are totally geodesic submanifolds.\qed
}\normalsize\\

Combining this theorem along with Corollary \ref{Coro-MaxIsoInteg} we find that if a maximally isotropic distribution is integrable then its leaves are totally geodesic, a known result that was proved in \cite{TaChabert-KY}. Particularly, if the signature of the manifold is Lorentzian and $I$ is a maximally isotropic distribution then $\dim(I\cap \overline{I})=1$, with $\overline{I}$ denoting the complex conjugate of the  distribution $I$ \cite{Simple Spinors}. Thus, if $I$ is integrable then $(I\cap \overline{I})$ generates a null geodesic congruence. In four dimensions, $2n= 4$, if the Ricci tensor vanishes then this geodesic congruence is shear-free and the Weyl tensor is algebraically special \cite{art2}.

%%%%%%%%%%%%%%%%%%%%
%%\begin{corollary}
%%If $I$ is a maximal isotropic distribution that is integrable then the leaves of this distribution are totally geodesic.
%%\end{corollary}
%%This corollary is an immediate consequence of theorem 3 and corollary 2.2, it is a known result that was proved before in \cite{TaChabert-KY}. If the %%signature of the manifold is Lorentzian and $I$ maximal then $\dim(I\cap \overline{I})=1$ \cite{Simple Spinors}, with $\overline{I}$ denoting the complex %%conjugate of the maximal isotropic distribution $I$. So, as observed in \cite{TaChabert-KY}, if $I$ is integrable in a manifold with Lorentzian signature %%then $(I\cap \overline{I})$ generates a null geodesic congruence\footnote{When $\dim(M)=4$ and the Ricci tensor vanishes this geodesic congruence is %%shear-free and the Weyl tensor is algebraically special \cite{art2}.}.
%%%%%%%%%%%%%%%%%%%%%%%%%%%%%

Note that all manipulations of the present section assumed just that the spinorial connection $\hat{\nabla}_a$ obeys to the Leibniz rule with respect to the Clifford action, meaning that Eq. (\ref{Leibniz-Cliff}) holds. Thus, in the above results the $1$-form $A$ of Eq. (\ref{Omega-Def}) is arbitrary. In particular, one can use this freedom and Corollary \ref{Coro-MaxIsoInteg} to prove that: If $\hat{\varphi}$ is a pure spinor that generates an integrable maximally isotropic distribution then, it is always possible to choose $A_a$ in (\ref{Omega-Def}) to be such that $\hat{\nabla}_X\,\hat{\varphi}= 0$ for all $X$ tangent to this distribution. Moreover, by means of Eq. (\ref{Curvature-Action}), one can verify that this required $1$-form can be a pure gauge, $A_a= \partial_a\lambda$ for some function $\lambda$, if, and only if, $\hat{\mathfrak{R}}_{XY}\,\hat{\varphi}= 0$ for all $X,Y$ tangent to the distribution.

%%%%%%%%%%%%%%%%%%%%%%%%%%%%%%%%%%%%%%%%%%%%%%%%%%%%%%%%%%%%%%%%%%%%%%%%%%%%%%%%%%%%%%%%%%%%%%%%%%%%%
%%%%%%%%%%%%%%%%%%%%%%%%%%%%%%%%%%%%%%%%%%%%%%%%%%%%%%%%%%%%%%%%%%%%%%%%%%%%%%%%%%%%%%%%%%%%%%%%%%%%%
%%%%%%%%%%%%%%%%%%%%%%%%%%%%%%%%%%%%%%%%%%%%%%%%%%%%%%%%%%%%%%%%%%%%%%%%%%%%%%%%%%%%%%%%%%%%%%%%%%%%%
%%%%%%%%%%%%%%%%%%%%%%%%%%%%%%%%%%%%%%%%%%%%%%%%%%%%%%%%%%%%%%%%%%%%%%%%%%%%%%%%%%%%%%%%%%%%%%%%%%%%%
%%%%%%%%%%%%%%%%%%%%%%%%%%%%%%%%%%%%%%%%%%%%%%%%%%%%%%%%%%%%%%%%%%%%%%%%%%%%%%%%%%%%%%%%%%%%%%%%%%%%%
%%%%%%%%%%%%%%%%%%%%%%%%%%%%%%%%%%%%%%%%%%%%%%%%%%%%%%%%%%%%%%%%%%%%%%%%%%%%%%%%%%%%%%%%%%%%%%%%%%%%%
%%%%%%%%%%%%%%%%%%%%%%%%%%%%%%%%%%%%%%%%%%%%%%%%%%%%%%%%%%%%%%%%%%%%%%%%%%%%%%%%%%%%%%%%%%%%%%%%%%%%%
%%%%%%%%%%%%%%%%%%%%%%%%%%%%%%%%%%%%%%%%%%%%%%%%%%%%%%%%%%%%%%%%%%%%%%%%%%%%%%%%%%%%%%%%%%%%%%%%%%%%%

\section{Twistor Equation and Integrability of Maximally Isotropic Distributions}\label{section-twistor}

It is well-known that in four dimensions a pure spinor obeying to the twistor equation generates an integrable distribution of isotropic planes. The aim of the present section is to investigate whether an analogous property holds in higher dimensions. Namely, the following questions are going to be answered: Does a pure spinor obeying to the twistor equation necessarily generate an integrable maximally isotropic distribution? What about the converse, does a pure spinor generating an integrable distribution obeys to the twistor equation when we judiciously use the freedom in the choice of the spinorial connection?

If $\hat{\varphi}\neq 0$ is an arbitrary pure spinor then one can always make a convenient choice of frame such that $\hat{\varphi}= \hat{1}$. In order to facilitate the calculations, let us assume that this choice was made, meaning that $e_i\cdot\hat{\varphi}= 0$ and $\theta^i\cdot\hat{\varphi}= \hat{\theta}^i$. Where, as previously defined, the indices $i,j,k,l$ run from $1$ to $n$ while the indices $a,b$ belong to $\{1,2, \ldots, 2n\}$, with $2n$ being the dimension of the manifold. Then, one can easily prove the following relations:
\begin{align}
 \nonumber (e_i\theta^k\theta^l)\cdot \hat{\varphi} = (\delta^k_i\,\hat{\theta}^l - \delta^l_i\,\hat{\theta}^k ) \quad\;;&\;\quad  (e_je_i\theta^k\theta^l)\cdot \hat{\varphi} = (\delta^k_i\,\delta^l_j - \delta^l_i\,\delta^k_j )\,\hat{1} \\
  \nonumber \hat{\nabla}_a \hat{\varphi} = ( A_a \,+\, \frac{1}{2}\,& \omega_{aj}^{\ph{aj}j} ) \, \hat{1} + \omega_{aij}\,\hat{\theta}^{ji}\,.
\end{align}
Where Eqs. (\ref{D-SpinorBasis}) and (\ref{Omega-Def}) were used. By means of the above relations one can prove, after some algebra, that the action of the twistor operator in $\hat{\varphi}$ is given by:
%%%%%%%%%%%%
%%\begin{align}
%%\nonumber \hat{\bl{T}}_j\, \hat{\varphi} \,=\,& \frac{1}{n}\left[ (n-1)(A_j+\frac{1}{2}\omega_{jk}^{\ph{jk}k} ) \,+\, \omega^{i}_{\ph{i}ij} \right] %%\hat{1} \,+ \\
%%+\, &\frac{1}{n}\left[ (n-1)\omega_{jkl} \,+\, 2\omega_{[lk]j} \right] \,\hat{\theta}^l\wedge\hat{\theta}^k
%%\end{align}
%%%%%
%%%\frac{1}{2}\, \hat{\bl{T}}^j\, \hat{\varphi} \,=\,
%%%%%%%%%%%%
\begin{align}
  \nonumber \hat{\bl{T}}_j\, \hat{\varphi} \,=\,& \frac{1}{n}\left[ (n-1)( A_j+\frac{1}{2}\,\omega_{jk}^{\ph{jk}k} ) \,+\, \omega^{i}_{\ph{i}ij} \right] \hat{1} \,+ \\
  +\, &\frac{1}{n}\left[ (n-1)\omega_{jkl} \,+\, 2\,\omega_{[lk]j} \right] \,\hat{\theta}^{lk} \label{Tjfi}
  \\
 \nonumber &\\
 \nonumber \hat{\bl{T}}_{j+n}\, \hat{\varphi} \,=&\,  \frac{1}{2}\lef A^j -
 \frac{1}{2} \omega^{ji}_{\ph{ji}i} \rig\, \hat{1} \,+\, \frac{1}{n} \, \omega_{[ikl]}\, \hat{\theta}^{jikl}  \,+\\
 +\, & \left[\, \frac{1}{2}\, \omega^{j}_{\ph{j}kl} \,+\,\frac{1}{n}\, \delta^j_{[k}A_{l]}  \,+\,\frac{1}{2n} \,\delta^j_{[k}\omega_{l]i}^{\ph{l]i}i} \,+\,  \frac{1}{n}\,\omega^{i}_{\ph{i}i[k}\delta^j_{l]} \,\right]\, \hat{\theta}^{lk} \label{Tj+nfi}\,.
\end{align}
Where in the above equations $V_{[a_1a_2\ldots a_p]}$ means the anti-symmetric part of $V_{a_1a_2\ldots a_p}$. For instance, $V_{[ab]}= \frac{1}{2!}(V_{ab}- V_{ba})$. If $\hat{\varphi}$ is a twistor then the right hand side of both equations above must vanish for all $j$, let us analyse these two conditions separately. First note that if the dimension is different from six, $n\neq 3$, then the right hand side of (\ref{Tjfi}) vanishes if, and only if,
\begin{equation}\label{Aj}
  A_j\,=\, \frac{1}{n-1}\, \omega^{i}_{\ph{i}ji} - \frac{1}{2}  \omega_{jk}^{\ph{jk}k} \quad\;\; \textrm{and} \quad \;\;\omega_{ijk}\,=\, 0\,.
\end{equation}
While in six dimensions $\hat{\bl{T}}_j\hat{\varphi}= 0$ if, and only if,
$$ A_j\,=\, \frac{1}{3-1}\, \omega^{i}_{\ph{i}ji} - \frac{1}{2}  \omega_{jk}^{\ph{jk}k} \quad\;\; \textrm{and} \quad \;\;\omega_{ijk}\,=\, \omega_{[ijk]}\,. $$
Since the maximally isotropic distribution associated to the pure spinor $\hat{\varphi}= \hat{1}$ is spanned by $\{e_1, e_2, \ldots, e_n\}$ and this distribution is integrable if, and only if, $\omega_{ijk}= 0$ then the above conditions lead us to the following theorem:
\begin{theorem}
If a pure spinor $\hat{\varphi}$ generates an inetgrable maximally isotropic distribution then one can always choose the 1-form $A_a$ so that $\hat{\bl{T}}_X \hat{\varphi}= 0$ for all vector fields $X$ tangent to such distribution. Conversely, if the dimension is different from six and $\hat{\bl{T}}_X \hat{\varphi}= 0$ for all $X$ tangent to the maximally isotropic distribution generated by $\hat{\varphi}$ then this distribution is integrable.
\end{theorem}

Now, assume that we choose the $1$-form $A$ so that the components $A_j$ are just as in Eq. (\ref{Aj}). Then, inserting this choice into (\ref{Tj+nfi}) leads us to the following equation:
\begin{align}
\nonumber \hat{\bl{T}}_{j+n}\, \hat{\varphi} \,=&\,  \frac{1}{2}\lef A^j -
 \frac{1}{2} \omega^{ji}_{\ph{ji}i} \rig\, \hat{1} \,+\,  \frac{1}{n} \, \omega_{[ikl]}\, \hat{\theta}^{jikl} \,+\\
 +&\,  \left[\, \frac{1}{2}\, \omega^{j}_{\ph{j}kl} \,+\,  \frac{1}{n-1}\,\omega^{i}_{\ph{i}i[k}\delta^j_{l]} \,\right]\, \hat{\theta}^{lk} \,.  \label{Tj+nfi2}
\end{align}
So, in order for the equation $\hat{\bl{T}}_{j+n}\hat{\varphi}= 0$ to hold the three terms on the right hand side of the above equation must vanish. The first term can always be made to vanish by a suitable choice of the $1$-form $A$, namely we must set
$$ 2A_{j+n} \,=\,  A^j \,=\,  \frac{1}{2} \omega^{ji}_{\ph{ji}i} \,. $$
Due to the anti-symmetry in its indices, the spinor $\hat{\theta}^{jikl}$ is zero in four and six dimensions. Therefore, in these cases the second term on the right hand side of Eq. (\ref{Tj+nfi2}) is automatically zero. While if the dimension is greater or equal to eight this term will vanish if, and only if, $\omega_{[ijk]}= 0$. But it is worth recalling that if we assume that $\hat{\bl{T}}_j\hat{\varphi}= 0$ then this condition is already satisfied. Finally, the third term on the right hand side of Eq. (\ref{Tj+nfi2}) vanish if, and only if,
\begin{equation}\label{conditions-3}
 \omega^i_{\ph{i}jk}\,=\,0 \;\;\;\forall\; i\neq j,k \,\quad \textrm{and} \quad \, \omega^i_{\ph{i}ik}\,=\, \omega^j_{\ph{j}jk} \;(\textrm{no sum in }i,j ) \;\;\;\forall\; i\neq k \neq j   \,.
\end{equation}
Where it is worth reemphasizing that in the second condition above no sum is assumed in the repeated indices $i$ and $j$. Since in four dimensions $i,j,k\in \{1,2\}$ it follows that the conditions of (\ref{conditions-3}) are identically satisfied in this case. While in higher dimensions these constraints are non-trivial. In order to give some sort of insight on the meaning of the constraints displayed in (\ref{conditions-3}), let us register that the integrability of the maximally isotropic distribution $\textrm{Span}\{e_i\}$ along with (\ref{conditions-3}) is tantamount to the following restrictions:
\begin{gather}
  g(\nabla_Xe_i-\nabla_Ye_j,e_k)  \,=\, 0  \quad\; \forall \;\, X,Y\in\,\Gamma(TM) \textrm{ such that}\label{Condition-geom}\\
 \nonumber   \; g(X,e_i)\,=\,g(Y,e_j) \;\textrm{ and }\;  g(X,e_k)\,=\,g(Y,e_k)\,=\,0 \,.
\end{gather}

The table below summarizes the joint analysis of Eqs. (\ref{Tjfi}) and (\ref{Tj+nfi}), displaying the necessary and sufficient conditions for the pure spinor $\hat{\varphi}$ to be a twistor. Since the condition $\omega_{ijk}= 0$ is equivalent to the integrability of the maximally isotropic distribution generated by $\hat{\varphi}= \hat{1}$, such table leads us to Theorem \ref{Theo-Twistor}, presented on the sequence.
\begin{table}[!htbp]
\begin{center}
\begin{tabular}{c|c|c|c}
  \hline
  \hline
  Dimension & $A_j$ &  $A_{j+n}$ & Constraints over $\omega_{abc}$\\
  \hline
  4 &\; $\frac{1}{2-1} \omega^{i}_{\ph{i}ji} - \frac{1}{2}  \omega_{jk}^{\ph{jk}k}$ \;&\;  $\frac{1}{4} \omega^{ji}_{\ph{ji}i}$ \;& $\omega_{ijk}= 0$\\
  6 & \;$\frac{1}{3-1}\omega^{i}_{\ph{i}ji} - \frac{1}{2}  \omega_{jk}^{\ph{jk}k}$     \;& \; $\frac{1}{4} \omega^{ji}_{\ph{ji}i}$ \;& $\omega_{ijk}= \omega_{[ijk]}$\,, Eq. (\ref{conditions-3})\\
  $2n\geq 8$ &\; $\frac{1}{n-1} \omega^{i}_{\ph{i}ji} - \frac{1}{2}  \omega_{jk}^{\ph{jk}k}$ \;&\;  $\frac{1}{4} \omega^{ji}_{\ph{ji}i}$ \;& $\omega_{ijk}= 0$\,, Eq. (\ref{conditions-3})\\
  \hline
  \hline
\end{tabular}\caption{\footnotesize{Depending on the dimension, the pure spinor $\hat{\varphi}= \hat{1}$ is a twistor if, and only if, the components of the $1$-form $A$ and the connection coefficients $\omega_{abc}$ are as displayed in this table. Recall that the constraints in the last row of the last column can be replaced by the single Eq. (\ref{Condition-geom}). }}\label{Tab.Twistor}
\end{center}
\end{table}
\normalsize

\begin{theorem}\label{Theo-Twistor}
If $\hat{\varphi}\neq 0$ is a pure spinor then the following results hold depending on the dimension of the manifold:\\
(1) $2n=4\rightarrow$ If $\hat{\varphi}$ is a twistor then the maximally isotropic distribution generated by $\hat{\varphi}$ is integrable. Conversely, if $\hat{\varphi}$ generates an integrable distribution then one can always choose $A_a$ so that $\hat{\varphi}$ obeys to the twistor equation with respect to the connection of Eq. (\ref{Omega-Def}). \\
(2) $2n=6\rightarrow$ The fact that $\hat{\varphi}$ obeys to the twistor equation does not imply that the distribution generated by  $\hat{\varphi}$ is integrable. Conversely, the integrability of the distribution generated by  $\hat{\varphi}$ does not guarantee that one can arrange $A_a$ in order to make  $\hat{\varphi}$  obey to the twistor equation.\\
(3) $2n\geq8\rightarrow$ If $\hat{\varphi}$ is a twistor then the maximally isotropic distribution generated by $\hat{\varphi}$ is integrable. On the other hand, the integrability of the distribution generated by $\hat{\varphi}$ does not imply that one can find $A_a$ such that $\hat{\varphi}$ obeys to the twistor equation.
\end{theorem}

Hopefully, the results presented in the last two sections will be valuable to some branches differential geometry. Specially, since the null directions play an important role on the study of holonomy \cite{Baum,Leistner-Galaev} it follows that some applications on this subject might appear. Although the study of Euclidean restricted holonomy is well established \cite{Berger}, for other signatures some questions remain open. In particular, recently some progress has been accomplished in the Lorentzian case \cite{Leistner-Galaev}. It is also worth mentioning the interplay between holonomy groups and the twistor equation \cite{Lischewski}. Concerning physics, the formalism introduced in the present article might be useful for string theory and supergravity \cite{M-waves} as well as in the study of geometrical properties of Black-Holes \cite{BGCC:ARQ}.

\section*{Acknowledgments}
I want to thank CNPq (Conselho Nacional de Desenvolvimento Cient\'{\i}fico e Tecnol\'{o}gico - Brazil) and CAPES (Coordena\c{c}\~{a}o de Aperfei\c{c}oamento de Pessoal de N\'{\i}vel Superior - Brazil) for the financial support. I also wish to thank the anonymous referee for helpful comments.

\end{document}